\documentclass[12pt]{article}
\usepackage{graphicx,amssymb,amsmath,amsthm,cite}

\def\Del{\Delta}

\newcommand{\la}{\langle}
\newcommand{\ra}{\rangle}

\newcommand{\Spann}{{\mbox{\rm{span}}}}
\newcommand{\til}{\tilde}
\newcommand{\vt}{\til{v}}

\newtheorem{remark}{Remark}

\def\be{\begin{equation}}
\def\ee{\end{equation}}
\def\ben{\begin{eqnarray}}
\def\een{\end{eqnarray}}

\def\D{\mathcal{D}}
\def\R{\mathbb{R}}
\def\N{\mathbb{N}}
\begin{document}

\title{Sparse image representation by discrete cosine$/$spline  
based dictionaries}
\author{James Bowley and Laura Rebollo-Neira\\
Mathematics, Aston University,UK}
\maketitle
\subsection*{Abstract}
Mixed dictionaries generated by cosine and B-spline functions 
are considered. It is shown that, by highly nonlinear approaches such 
as Orthogonal Matching Pursuit, the discrete version of 
the proposed dictionaries yields a 
significant gain in the sparsity of an image representation.
\section{Introduction} 
Sparse representation of information is a central aim of 
data processing techniques. An usual first step of image 
processing applications, for instance, is 
to map the image onto a transformed space allowing for the reduction 
of the number of data points representing the image. Currently the 
most broadly used transforms for performing that task are the 
Discrete Cosine Transform (DCT) and Discrete Wavelet Transforms (DWT). 
An important reason for the popularity of both these transforms 
is the viability of their fast implementation. However, since parallel 
processing is becoming more powerful and accessible, alternative 
approaches for signal representation are being given increasing 
consideration. 
Emerging techniques address the matter in the following way: 
Given a signal $f\in \R^N$ find the decomposition 
$f=\sum_{i=1}^M c_{i} v_{\ell_i}$, where vectors 
$v_{\ell_i} \in \R^N ,\,i=1,\ldots,M$, usually called atoms, 
are a subset of a redundant set called a dictionary. 
Approximations of this type are highly nonlinear and are said to 
yield a sparse representation of the signal $f$ in terms of $M$ atoms 
if $M$ is considerably smaller than $N$. Available methodologies
for nonlinear approximations are known as Pursuit Strategies. This 
comprises Bases Pursuit \cite{CDS98} and Matching Pursuit 
like greedy algorithms \cite{MZ93}, including 
Orthogonal Matching Pursuit (OMP) and variations of it
\cite{PRK93,RL02,AR06,AR06a}. Another concern inherent to
highly non linear approximations is the design
of suitable dictionaries for representing some class of signals. 
Dictionaries arising by merging orthogonal
bases are theoretically studied in \cite{DH01,GN03}.
From a different perspective, approaches for learning dictionaries 
from large data sets are considered in \cite{OF97,AEB06}.
In this communication we present an alternative construction 
of dictionaries for representing natural images.
The proposed dictionaries are a mixture of discrete cosine 
and spline based dictionaries. We have found
by a good number of examples, some of which are 
presented here, that the resulting dictionary renders 
a considerable gain in sparsity in comparison 
to fast transforms such as DCT and DWT, at acceptable 
visual level (PSRN 40 dB). 

The paper is organized as follows: In Sec.~\ref{bspline}
we introduce the discrete B-spline based dictionaries which 
together with the discrete cosine ones form the large mixed 
dictionary we are proposing. In Sec.~\ref{omp} 
we discuss the implementation of the OMP approach that we have
used. The details of the actual process for dealing with 
images are given in Sec. \ref{exa} where results illustrating 
the capability of the proposed dictionaries to yield 
sparse representations by nonlinear approaches are provided. 

\section{B-spline based dictionaries}
\label{bspline}
The discrete dictionaries we discuss here are inspired by a 
general result holding for continuous spline 
spaces. Namely, that {\em{spline spaces on a closed interval
can be spanned by dictionaries of B-splines
of broader support than the corresponding B-spline basis functions 
\cite{AR05}}}. 

A {\em{partition}} of an interval $[c,d]$ is a finite set of points
$\Del:=\{x_i\}_{i=0}^{N+1}, N\in\N$ such that 
$c=x_0<x_1<\cdots<x_{N}<x_{N+1}=d$, which generates 
$N$ subintervals 
$I_i=[x_i,x_{i+1}), i=0,\dots,N-1$ and $I_N=[x_N,x_{N+1}]$. 
Representing by $\Pi_{m}$ the   
space of polynomials of
degree smaller than or equal to $m\in\N_0=\N\cup\{0\}$ and as
$C^{m}$ the space of functions having continuous derivatives up
to order $m$  (with $C^{0}$ the space of continuous functions) 
the spline space of order $m\ge 2$ on $[c,d]$, with single knots at the 
partition points, is define as 
$$S_m(\Del)=\{f\in C^{m-2}[c,d]\ : \ f|_{I_i}\in\Pi_{m-1},
i=0,\dots,N\},$$
where $f|_{I_i}$ indicates the
restriction of the function $f$ to the
interval ${I_i}$.   
In the case of equally spaced knots the corresponding splines are called 
cardinal. Moreover all the cardinal
B-splines of order $m$ can be obtained from
one cardinal B-spline $B(x)$ associated with the uniform simple 
knot sequence $0,1,\dots,m$.
Such a function is given as
\begin{equation}
B_m(x)=\frac{1}{m!}\sum_{i=0}^m(-1)^i\binom{m}{i}(x-i)^{m-1}_+,
\end{equation}
where $(x-i)^{m-1}_+$ is equal to $(x-i)^{m-1}$ if $x-i>0$ and 0 otherwise.
\begin{figure}
 \begin{center}
\includegraphics[width=8cm]{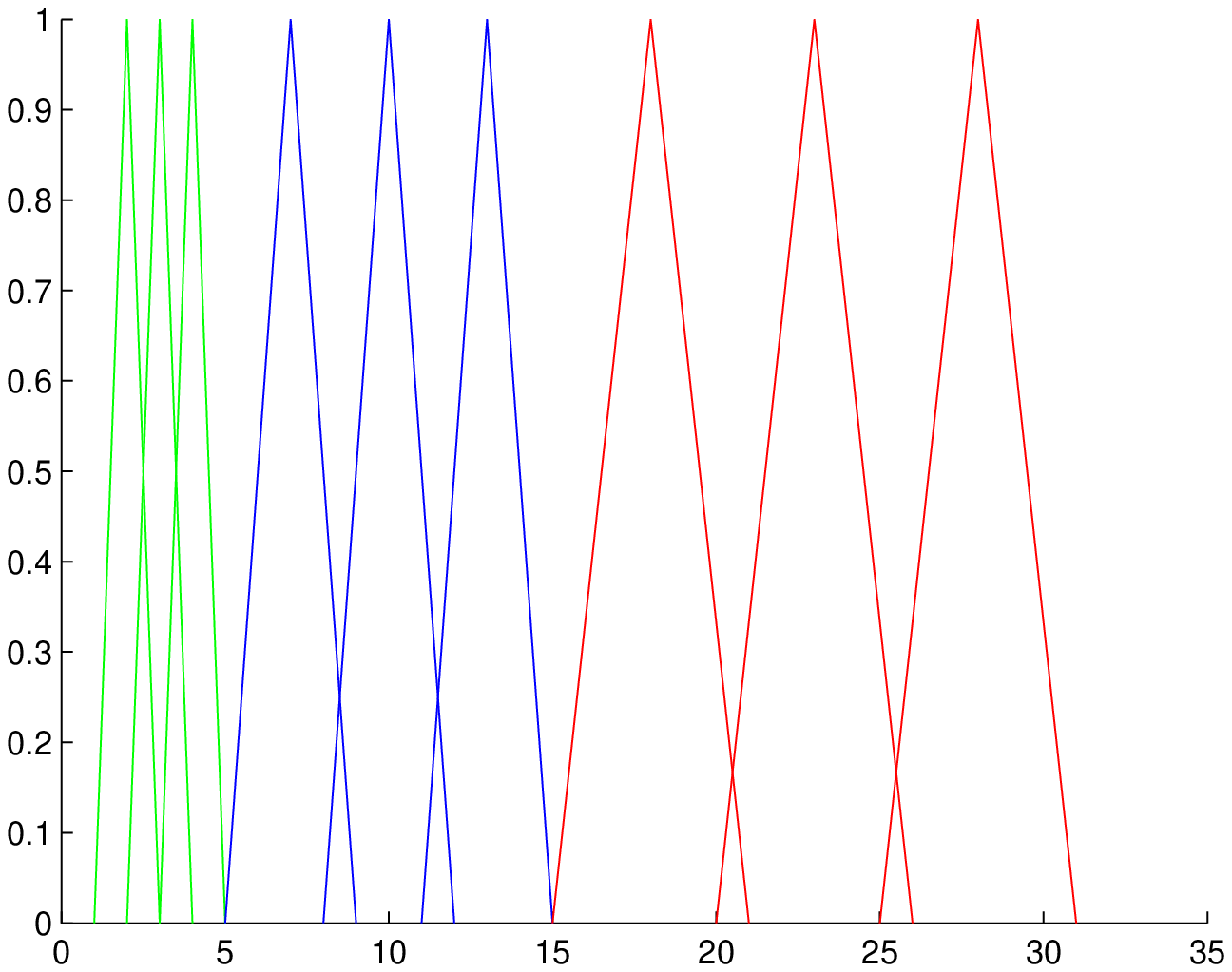}
\includegraphics[width=8cm]{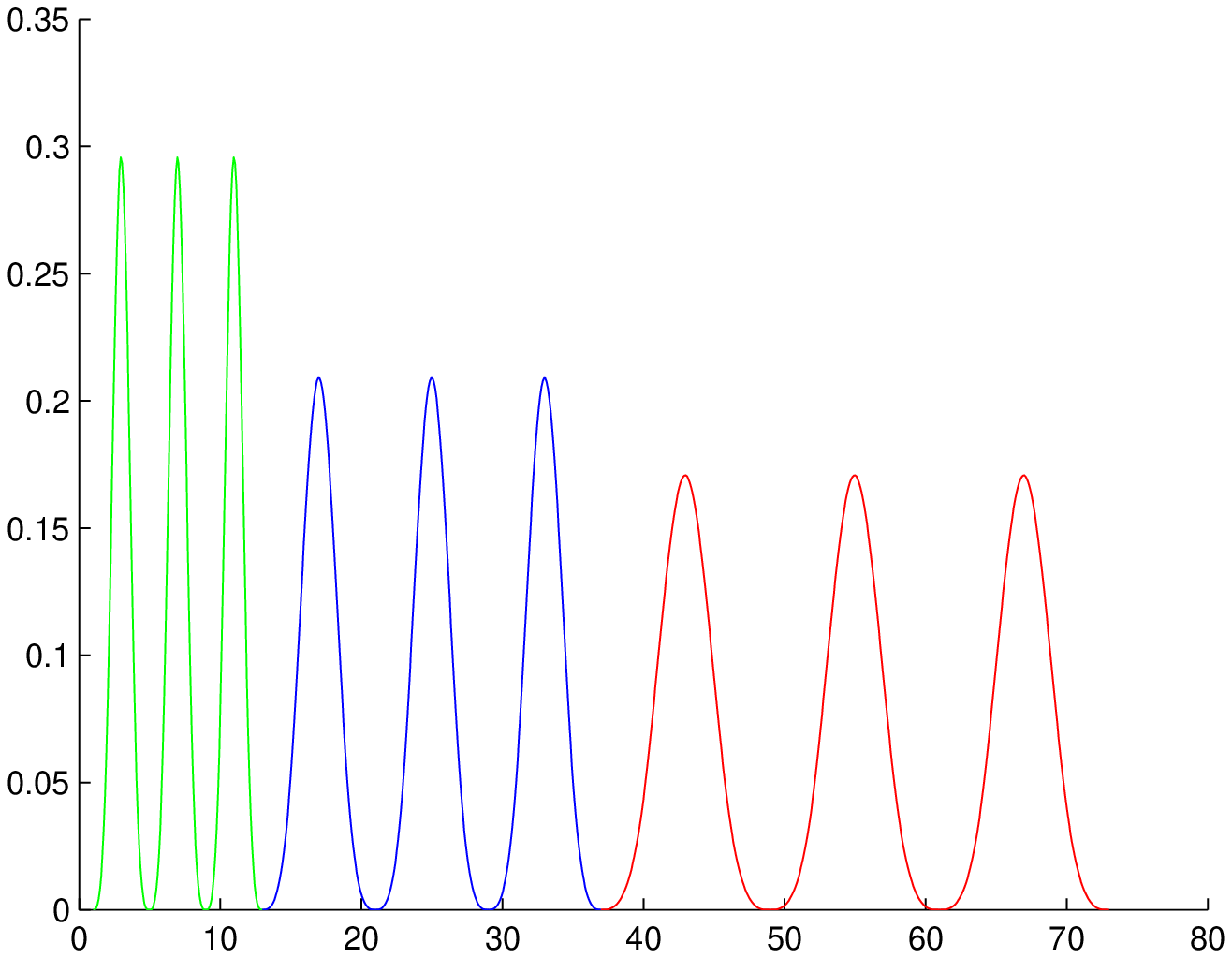}
\end{center}
\caption{Splines taken from
dictionaries spanning the same space. Linear splines (left)
Cubic splines (right)}
\label{splines}
\end{figure}
We shall focus on the particular cases corresponding to 
$m=2$ and $m=4$. For $m=2$ the cardinal spline space $S_2(\Del)$ is the 
space of piece wise linear functions and can be spanned by
a linear B-spline basis arising by translating the prototype 
function known as `hat' function. The first 3 functions 
in the left graph of Fig.~\ref{splines} are 3 consecutive 
linear B-spline basis functions. The 3 middle functions in 
the same graph are linear B-spline functions of broader support 
taking from a dictionary spanning the same space as the basis. 
The last 3 functions are taking from another dictionary 
for the same space. Details on how to build these 
dictionaries are given in \cite{AR05}.
The basis and dictionary functions equivalent to the 
ones in the left graph of Fig.~\ref{splines}, but
for cubic spline spaces corresponding to $m=4$, 
are given in the right graph of the same Figure. 

\begin{figure}
 \begin{center}
\includegraphics[width=8cm]{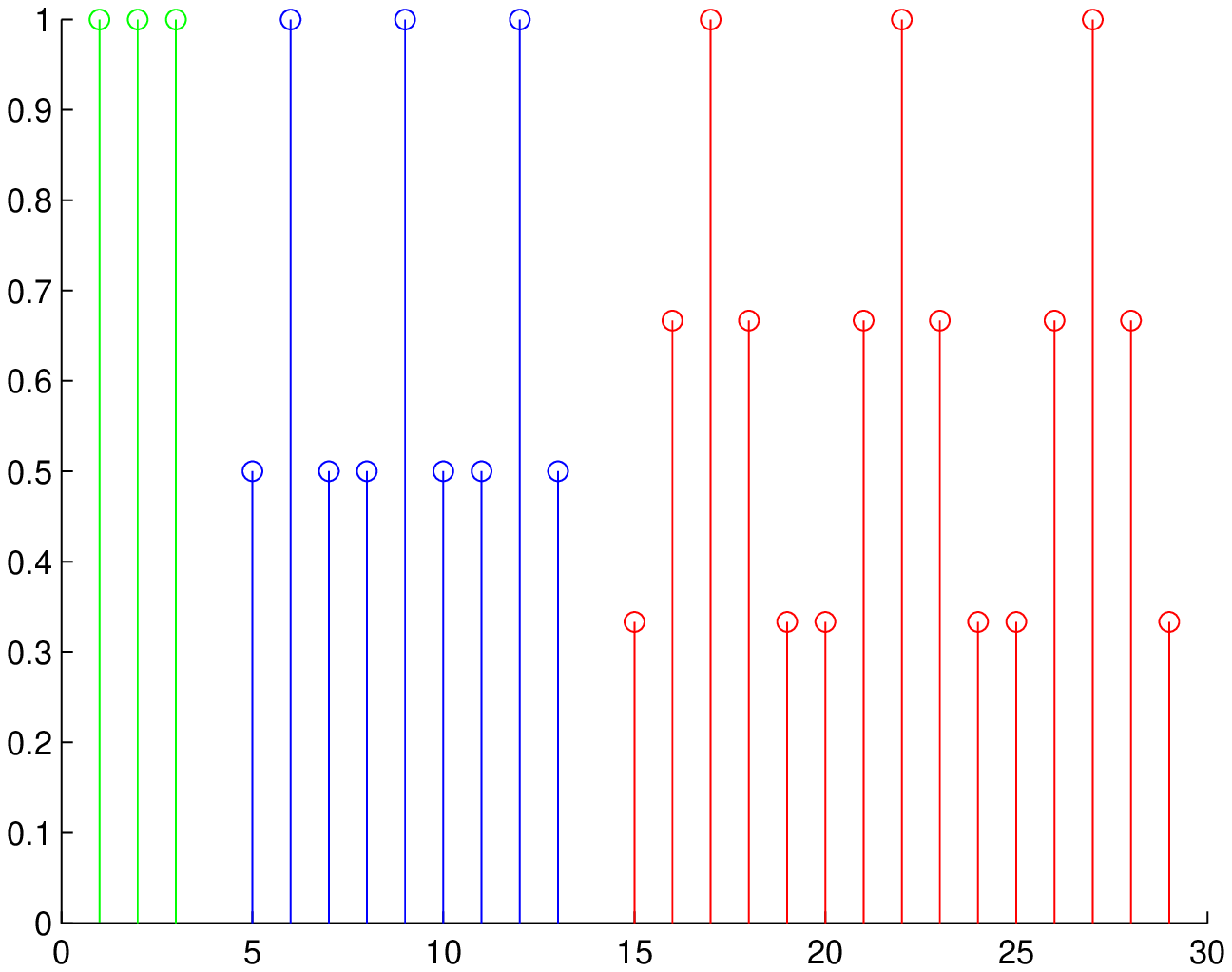}
\includegraphics[width=8cm]{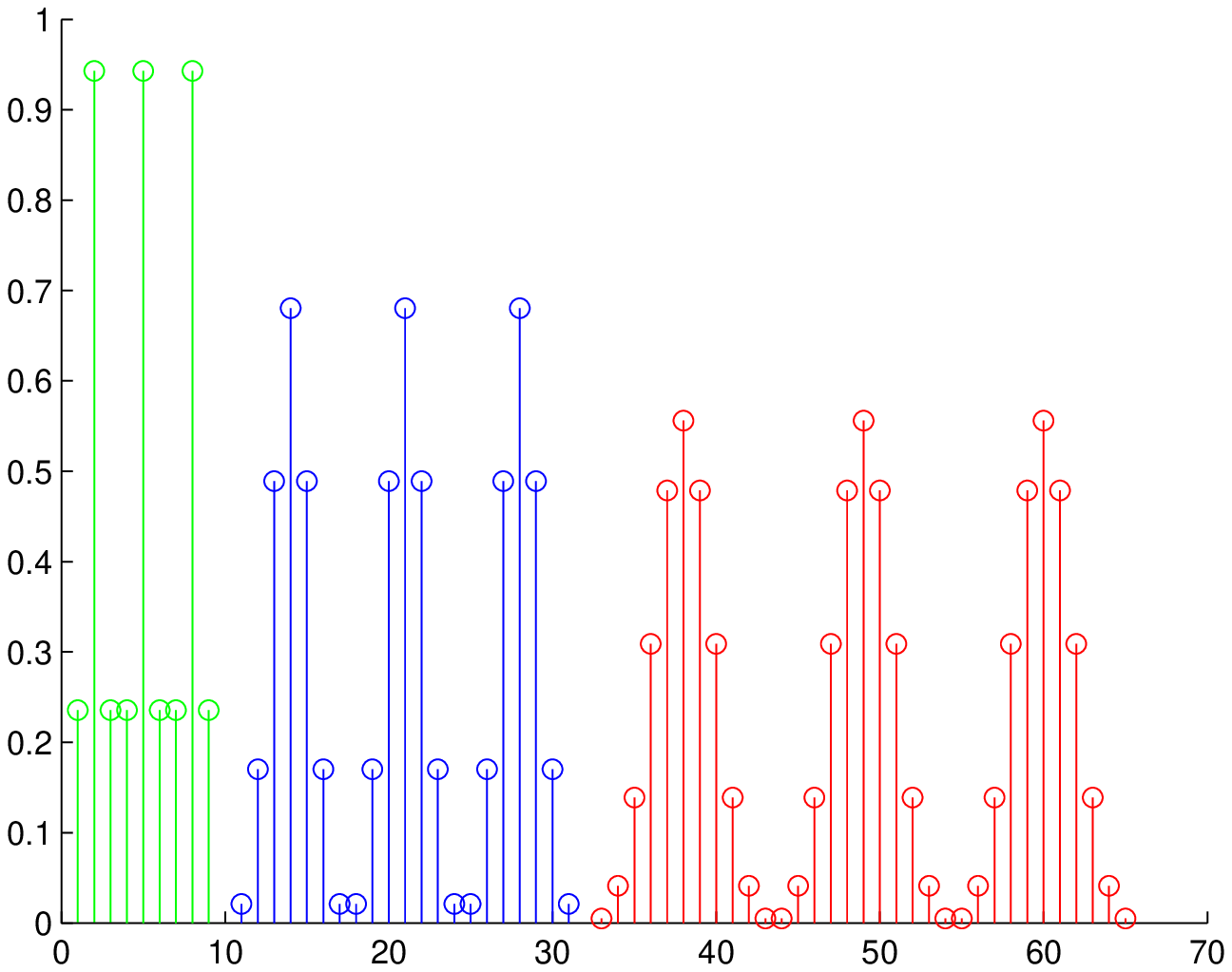}
\end{center}
\caption{Tree splines taken from
dictionaries spanning the same space. Linear splines (left)
Cubic splines (right)}
\label{splinesa}
\end{figure}

For constructing redundant dictionaries suitable for 
processing images by nonlinear techniques we need to make sure that
the dictionaries can be processed with digital computers 
having the existing memory capacity. Thus, we need to
a)discretise the functions to obtain adequate Euclidean 
vectors and b)restrict the functions to small intervals allowing
for processing the images in small blocks.
We carry out the discretization by taking the value of a
prototype function only at the knots (cf. small circles 
in graphs Fig. ~\ref{splinesa}) and translating that 
prototype one sampling point at each translation step. In regard 
to the boundaries one may take different routes:
A possibility is to adopt periodicity (cyclic boundary conditions)
and other apply the `cut off' approach and keep all the vectors 
whose support has nonzero intersection with the interval being considered. 
The former would leave a basis for the corresponding Euclidean 
space and the later a redundant dictionary.
\begin{remark}
We notice that by the proposed discretization 
the hat B-spline basis for the corresponding 
interval becomes the standard 
Euclidean basis for either boundary conditions. 
By discretizing the hats of broader support 
the samples preserve the hat shape.
\end{remark}
Obviously for a finite dimension Euclidean space we can construct 
arbitrary dictionaries. In particular,
different B-spline based dictionaries each of which comprising 
vectors of different support. 
Furthermore, we can include vectors of different support by 
merging dictionaries. There is, of course, a compromise between 
redundancy and complexity that needs to be considered. The discussion 
of such a tradeoff is postponed to Sec. \ref{exa}, where
the numerical examples are described.

From discrete unidimensional B-spline based dictionaries we 
obtain bidimensional ones simply by taking tensor product. 
Actually in Sec. \ref{exa} we consider a dictionary consisting 
of unidimensional cosine and B-spline based vectors, 
of redundancy approximately five, and build the bidimensional 
one by taking the tensor product of the whole dictionary with 
itself. 
\section{Implementation of the greedy algorithm OMP}
\label{omp}
The OMP technique \cite{PRK93} is an adaptive greedy strategy for 
selecting atoms which evolves as follows:
Let $f\in \R^N$ be a given signal and $\{v_i\}_{i=1}^L$ a given redundant 
dictionary. Setting $R^1=f$ at iteration $k+1$ the OMP algorithm
selects the atom, $v_{l_{k+1}}$ say, as the one minimizing the 
absolute values  of the inner products
$\la v_i ,R^{k} \ra,\, i=1,\ldots,L$, i.e.,
\be
\label{sel}
v_{l_{k+1}}= \arg\max\limits_{i \in J}|\la v_i ,R^{k} \ra |,\quad
\text{where}\quad R^{k}= f - \sum_{i=1}^{k} c_i^k v_{\ell_i},
\ee
and $J$ is the set of indices labeling the dictionary's atoms.
The coefficients $c_i^k,\,i=1,\ldots,k$ in the above decomposition 
are such that 
$||f- R^{k}||^2$ is minimum, which is equivalent to requesting
$R^{k}= \hat{P}_{V_k} f$, where $\hat{P}_{V_k}$ is the orthogonal 
projection operator onto $V_k= \Spann{\{v_{\ell_i}\}_{i=1}^k}$.
We base our implementation for determining the 
coefficients $c_i,\,i=1,\ldots,k$ on Gram Schmidt 
orthogonalization with re-orthogonalization, 
and recursive biorthogonalization.
Basically, at each iteration we update  the vectors
$$\vt_i^{k+1}= \vt_i^k - \vt_{k+1}^{k+1}\la v_{\ell_{k+1}}, \vt_i^k \ra,$$
where $\vt_{k+1}^{k+1}= q_{k+1}/\| q_{k+1}\|^2$, with 
$q_{k+1}= v_{\ell_{k+1}} - \hat{P}_{V_k} q_{k+1}$ and $q_1=v_{\ell_1}$. 
One reorthogonalization step implies to recalculate 
$q_{k+1}$ as $q_{k+1}= q_{k+1}-\hat{P}_{V_k} q_{k+1}$. The projector
$\hat{P}_{V_k}$ is here computed as $\hat{P}_{V_k}=Q_k Q_k^\ast$ where the 
$k$-columns of matrix $Q_k$ are the vectors $q_{i}/\| q_{i}\|,\,
i=1,\ldots,k$ and $Q_k^\ast$ indicates the transpose conjugate of 
$Q_k$. However, to calculate the coefficients of the 
linear superposition we express the projectors as
$\hat{P}_{V_k}={A_k B_k^\ast}$ where the $k$-columns of matrix 
$A_k$ are the selected vectors and the $k$-columns of matrix $B_k$ 
are the vectors $\vt_i^k,\,i=1,\ldots,k$. Thus, the required 
coefficients arise from the inner products
$c_i^k= \la \vt_i^k, f\ra,\, i=1,\ldots,k$. 
Details on this type of implementation are given in 
\cite{RL02,AR06} and the code can be found at
\cite{Webpage2}. Moreover, as will be discussed in the next 
section, the fact that we deal with dictionaries 
involving DC and supported atoms reduces the general 
complexity of the OMP method. 

\section{Sparse image representation by  Discrete Cosine and B-spline 
based dictionaries}
\label{exa}
Here we present the examples of the gain in sparsity that 
are achieved by using dictionaries which are the union of 
Discrete Cosine an B-Spline based dictionaries.  
In order to apply the OMP approach on an image using dictionaries 
of theses types is necessary to divide the 
image into blocks. The size of all the test images we consider is 
$512 \times 512$ and we divide each image into blocks of 
$16\times16$ pixels. For approximating each block we 
first construct the dictionaries
${\D_1}$ and ${\D_i},\,i=2,\ldots,4$  
defined as follows
\begin{itemize}
\item Discrete Cosine Dictionary

$$\D_1=\{c_i\cos(\frac{\pi(2j-1)(i-1)}{4L}),
\,j=1,\ldots,N \}_{i=1}^M,$$ 
with $c_i=\frac{1}{\sqrt{L}}\,i=1,\ldots,M$ and 
$c_1=\frac{1}{\sqrt{2L}}$ a normalization factor. 

\item Discrete B-Spline based dictionaries

$$\D_k=\{w_iB^k_m(j-i)|L; j=1,\ldots,L \}_{i=1}^{M_k},$$

where the notation $B_m(j-i)|L$ indicates the restriction to be 
an array of size $L$, indices $k=2,\ldots,4$  
label the dictionaries of different support 
and $w_i,\,i=1,\ldots,M_k$, with $M_k$ equal to the 
number of atoms in dictionary $k$, are normalization constants. 
Considerations are limited to the cases  $m=2$ 
(hat atoms) and $m=4$ (atoms arising by discretizing cubic 
splines). Because we adopt the cut off approach for the boundary, 
the numbers $M_k$ of total atoms in the $k$th-dictionary
varies according to the atom's support. For the linear
spline based dictionaries the corresponding supports 
are: 1, 3, and 5, while for the cubic are 3, 7, and 11. 
\end{itemize}
With these dictionaries we construct the tensor product dictionary 
$\D= \D_i\otimes \D_j,\,i,j=i,\ldots,4$ which implies, 
approximately, a redundancy of five. However, 
the redundancy does not modify significantly the complexity order. The 
complexity of applying the OMP approach is dominated by the 
evaluation, at each iteration, of the inner products 
between the residual and the dictionary atoms. In the 
case of the proposed dictionaries the inner product with the 
DC dictionaries can be implemented by fast DCT and 
the complexity in computing the inner products with the
other atoms depends on the atoms support (cf \eqref{sel}). 
By denoting $d_i$ to the support of the spline based dictionary $i$, the 
complexity of computing the inner products   
at the selection step \eqref{sel} is
$$O((2N)^2 \log_2 2N) + O(N \sum_{i=2}^4 d_i M_i).$$
For the supported atoms we are considering the number of 
atoms $M_i,\,i=2,\dots,4$ is not much larger than $N$. 
Moreover, for the hats dictionaries $d_1=1$, $d_2=3$ $d_3=5$, so that
the redundancy do not affect so much the stepwise complexity 
dominated by the inner products of the residual with all the 
dictionary's atoms. In addition  which can be also accelerated by parallel 
calculations.
Of course, the total complexity depends on the sparsity, as 
the complexity for selecting each atom has to be multiplied by
the number of selected atoms, which is the feature of stepwise
Pursuit Strategies. Nevertheless, the fact that the image is processed
in small blocks leaves room for fast implementation by parallel 
processing.

As can be observed in Table I, the  performance in sparsity that 
is achieved  with the proposed dictionaries is slightly better 
when using hats dictionaries, but with both dictionaries the sparsity 
in representing the six test images is about double than that 
yielded by faster nonlinear techniques such as DCT and WT. 

\begin{figure}
\begin{center}
\includegraphics[width=4cm]{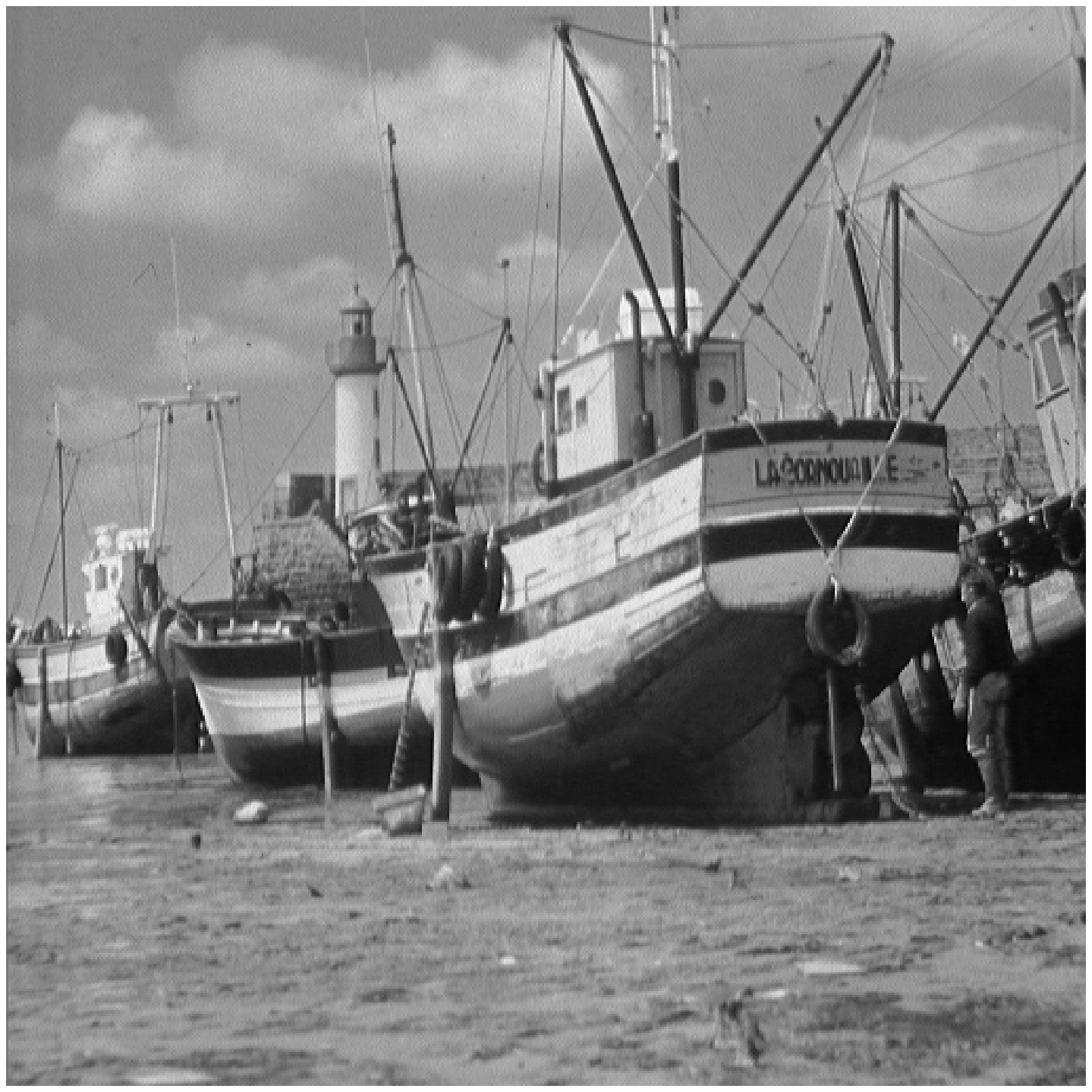}\hspace{-1cm}
\includegraphics[width=4cm]{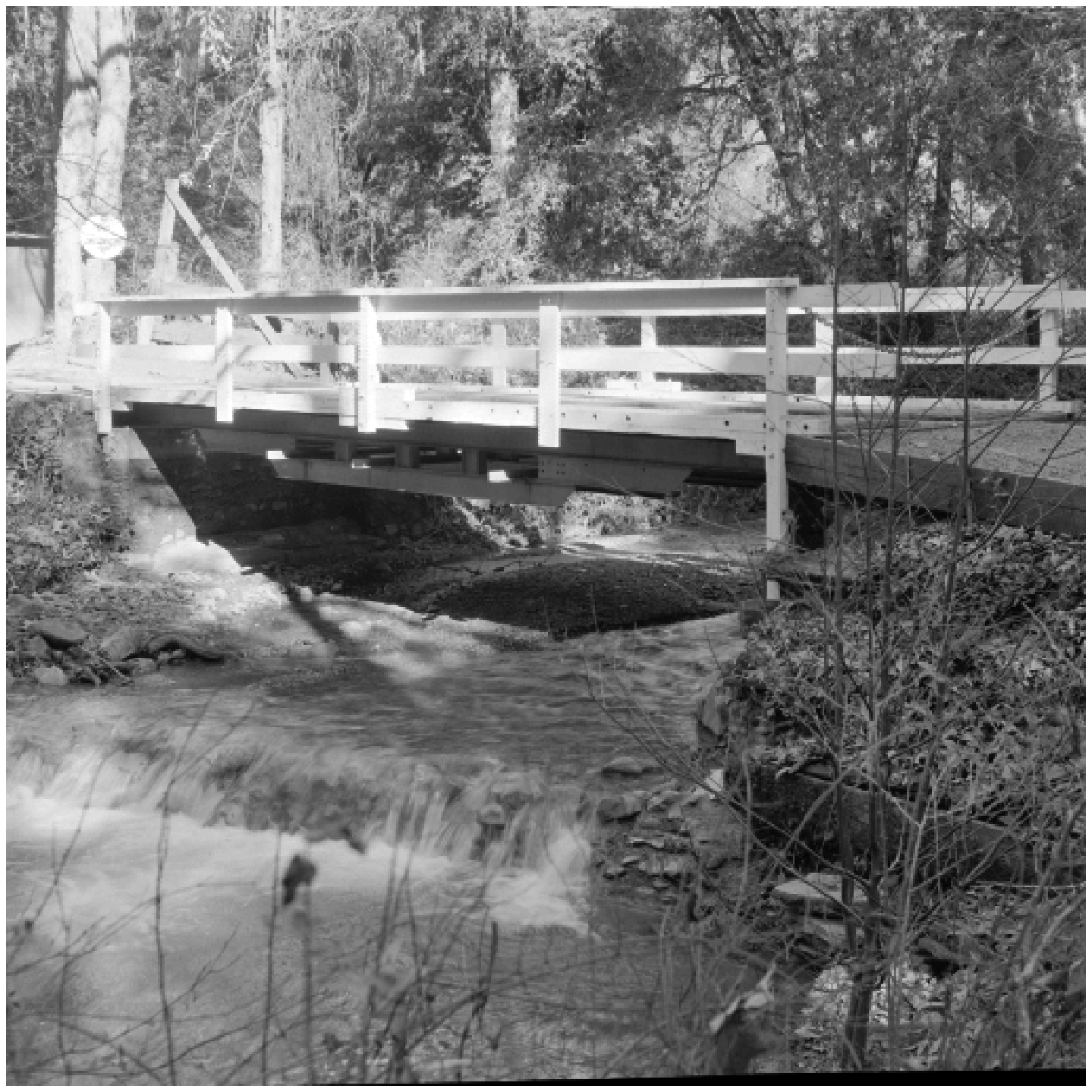}\\
\includegraphics[width=4cm]{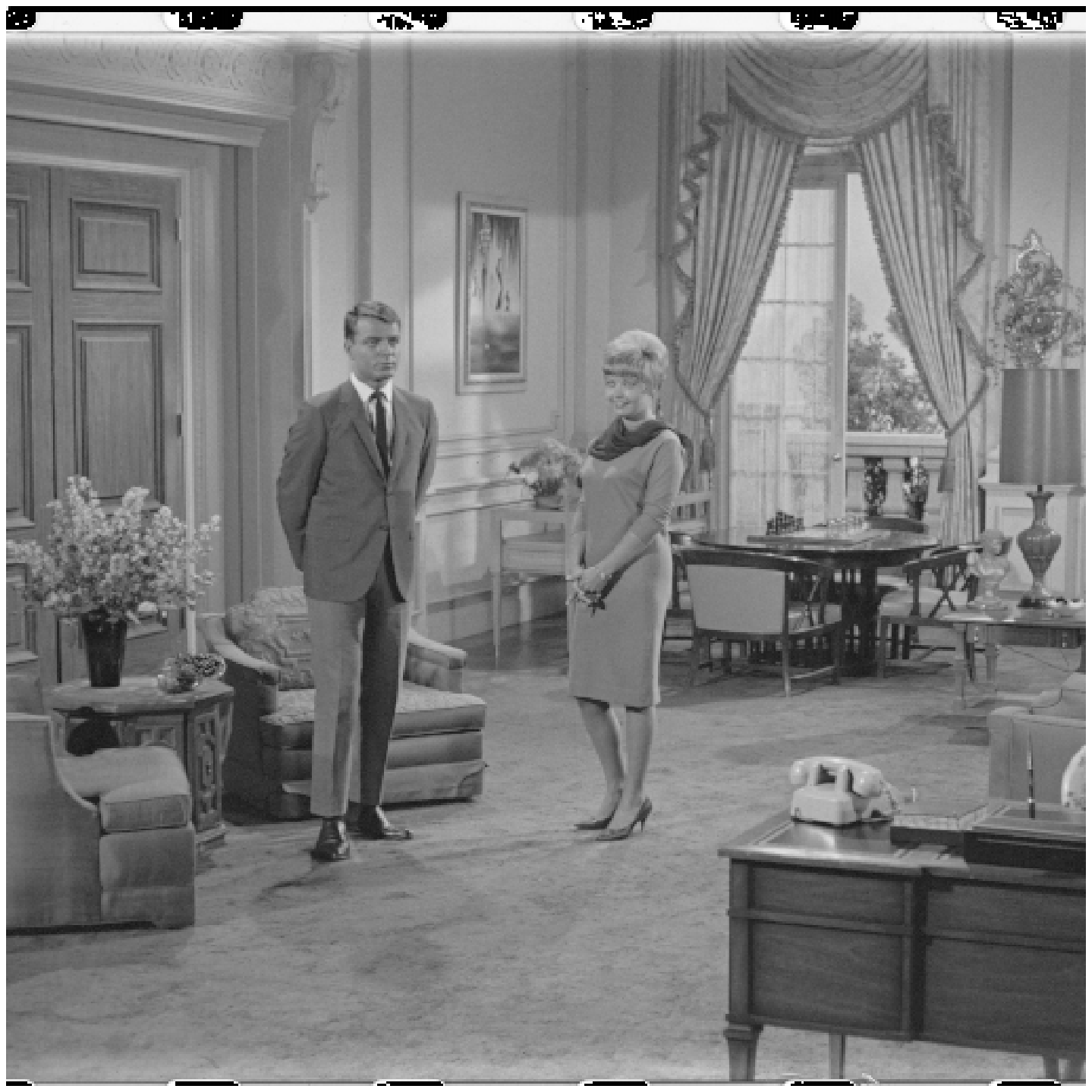}\hspace{-1cm}
\includegraphics[width=4cm]{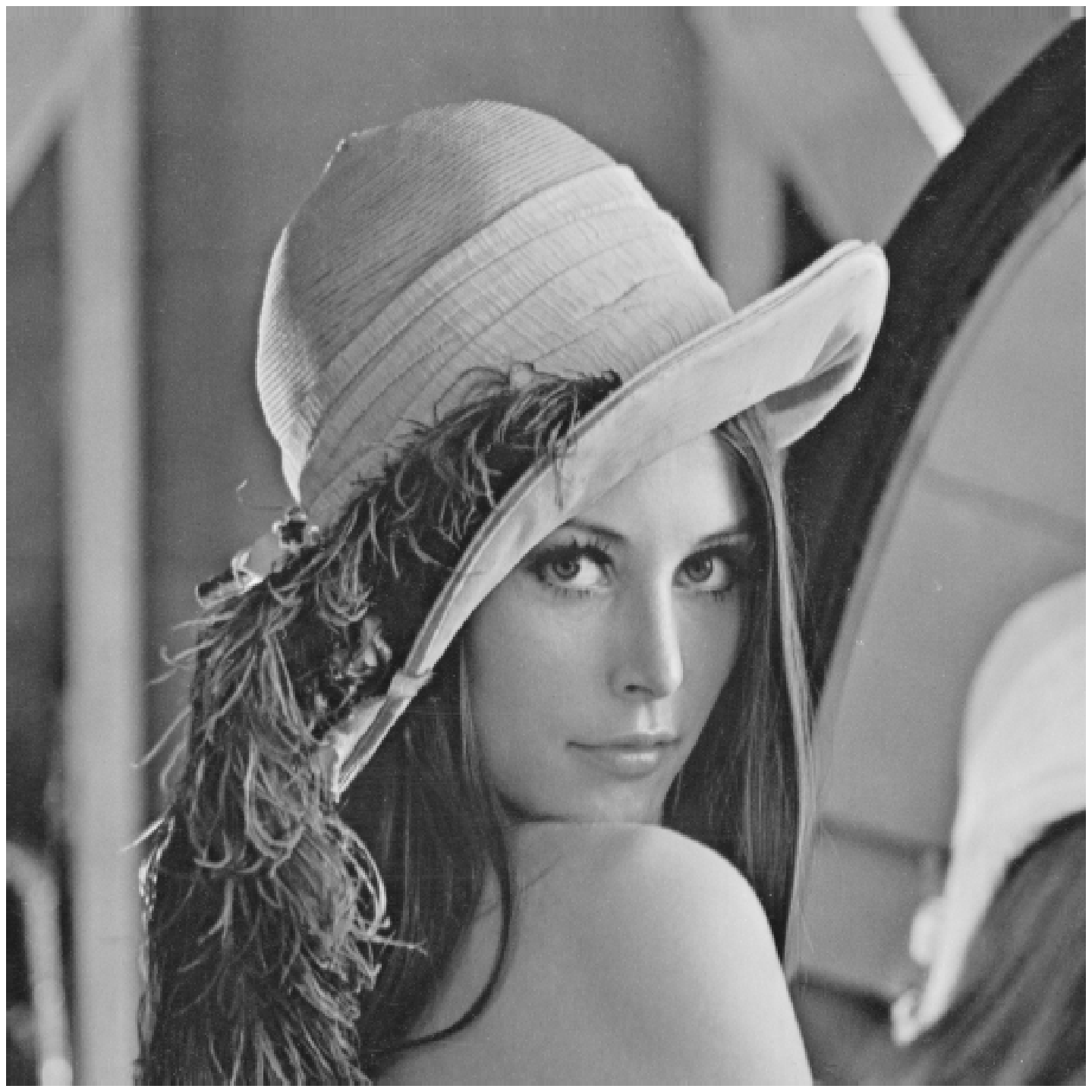}\\
\includegraphics[width=4cm]{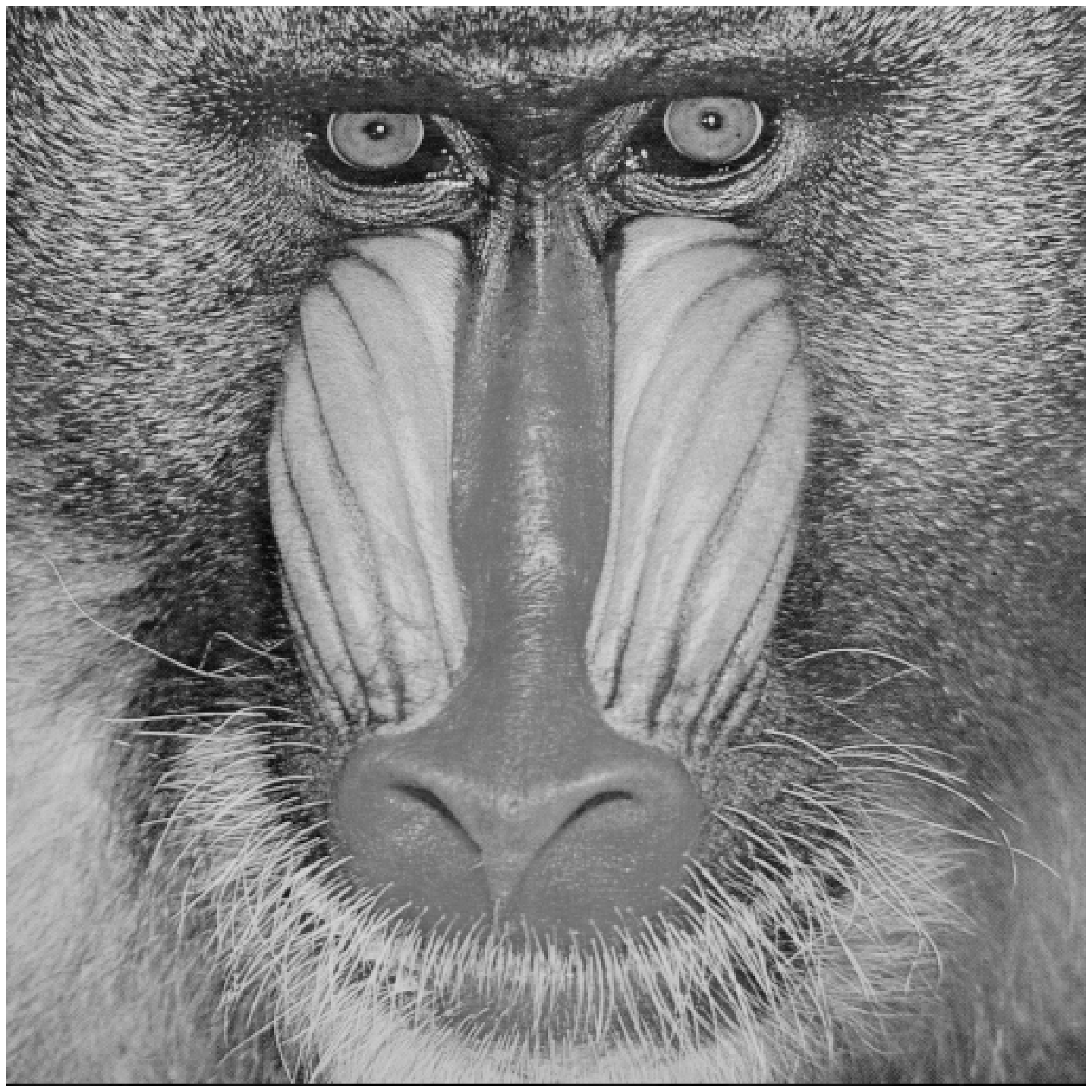}\hspace{-1cm}
\includegraphics[width=4cm]{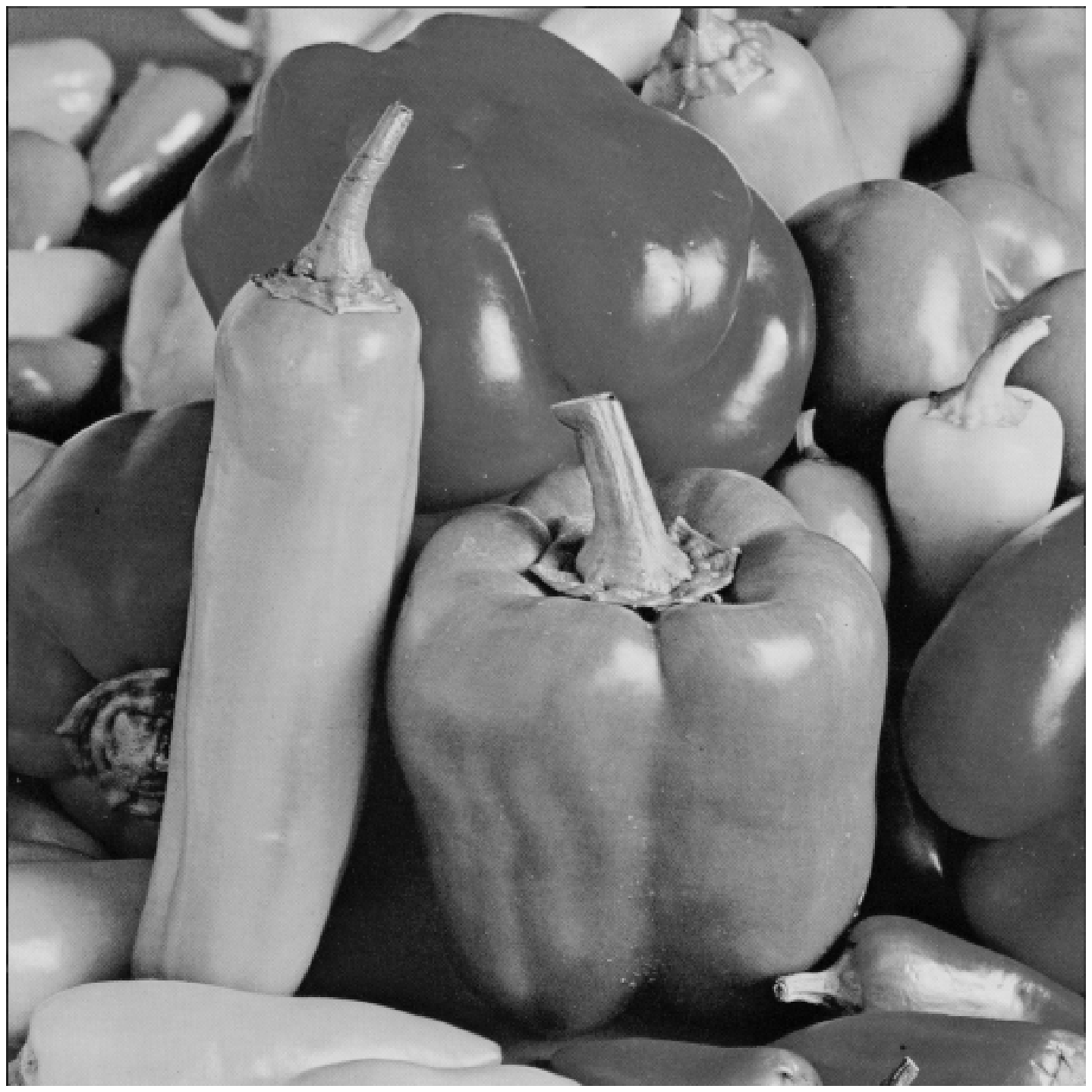}
 \end{center}
 \caption{The six test images from left to right,
top to bottom: Boat, Bridge, Film clip, Lena, Mandril, Peppers}
\label{test_images}
\end{figure}
\begin{table}
\begin{center}
\begin{tabular}{ | l || c | c | c | c | c | }
\hline
Image & DCT2 $\cup$ Linear Splines & DCT2 $\cup$ Cubic Splines & DCT &
Wavelets\\ \hline \hline
Boat & 7.05 & 6.89  & 3.63 & 3.65 \\ \hline
Bridge & 4.24  & 3.97 & 2.06 & 2.2 \\ \hline
Film & 9.72 & 9.26 & 4.53 & 4.8 \\ \hline
Lena & 11.78 & 11.7 & 6.5 & 6.97 \\ \hline
Mandril & 3.72 & 3.5 & 1.91 & 1.90 \\ \hline
Peppers & 8.9 & 8.62 & 4.36 & 3.39 \\ \hline
\end{tabular}
\caption{Compression ratio (corresponding to PNSR=40 dB)
achieved by each dictionary. The first column
corresponds to the dictionary $\D_1$ composed  of
DC redundancy 2 and linear spline atoms of support 1, 3 and 5.
The second column corresponds to dictionary  $\D_2$
and cubic spline atoms of support 3, 7 and 11. The
third column corresponds to the result obtained by
nonlinear selection of DCT coefficients. The last column
is the compression ratio produced by the
Cohen-Daubechies-Feauveau 9/7 wavelet transform computed
with the software WaveletCDF97
by thresholding coefficients so as to achieve the
required PNSR or 40 dB)}
\end{center}
\end{table}
\section{Conclusions}
\label{con}
Mixed DC and spline based dictionaries for sparse 
image representation have been introduced. It was shown that, 
comparing with fast nonlinear DCT and WT approaches, the proposed
dictionaries yield a significant gain in sparsity.
The complexity analysis and the fact that the processing 
is suitable for parallel computing lead to conclude that 
the proposed dictionaries can be of assistance to those image
processing applications that benefit from the sparsity 
property of a representations.  
\bibliographystyle{elsart-num}
\bibliography{revbib}
\end{document}